\documentclass{amsart}
\usepackage[all]{xy}
\usepackage{color}
\usepackage{amsthm}
\usepackage{amssymb}
\usepackage[colorlinks=true]{hyperref}
\numberwithin{equation}{section}

\newtheorem{theorem}[equation]{Theorem}
\newtheorem*{theorem*}{Theorem}

\newtheorem*{conjecture*}{Mamma Conjecture}
\newtheorem*{conjecture1*}{Mamma Conjecture (revisited)}
\newtheorem{proposition}[equation]{Proposition}

\newtheorem*{corollary*}{Corollary}

\theoremstyle{remark}
\newtheorem{definition}[equation]{Definition}

\newtheorem{example}[equation]{Example}

\theoremstyle{remark}
\newtheorem{remark}[equation]{Remark}

\newcommand{\cA}{{\mathcal A}}
\newcommand{\cB}{{\mathcal B}}
\newcommand{\cC}{{\mathcal C}}
\newcommand{\cD}{{\mathcal D}}

\newcommand{\cJ}{{\mathcal J}}

\newcommand{\cM}{{\mathcal M}}
\newcommand{\cN}{{\mathcal N}}
\newcommand{\cO}{{\mathcal O}}

\newcommand{\cT}{{\mathcal T}}
\newcommand{\cU}{{\mathcal U}}

\newcommand{\Spt}{\mathsf{Spt}}% Spectra
\newcommand{\Nat}{\mathsf{Nat}} % natural transformations

\newcommand{\bbA}{\mathbb{A}}

\newcommand{\bbD}{\mathbb{D}}

\newcommand{\bbL}{\mathbb{L}}

\newcommand{\bbP}{\mathbb{P}}
\newcommand{\bbR}{\mathbb{R}}

\newcommand{\bbQ}{\mathbb{Q}}
\newcommand{\bbZ}{\mathbb{Z}}

\DeclareMathOperator{\Id}{Id}

\DeclareMathOperator{\Mot}{Mot}
\DeclareMathOperator{\NChow}{NChow} % category of noncommutative Chow motives
\DeclareMathOperator{\NNum}{NNum} % category of noncommutative numerical motives
\DeclareMathOperator{\Chow}{Chow} % category of Chow motives
\DeclareMathOperator{\Num}{Num} % category of numerical motives

\DeclareMathOperator{\Mix}{KMM} % Kontsevich's category of noncommutative motives
\DeclareMathOperator{\Fun}{Fun} % Functor category
\newcommand{\bbK}{I\mspace{-6.mu}K}

\newcommand{\dgcat}{\mathsf{dgcat}}

\newcommand{\perf}{\mathsf{perf}}

\newcommand{\dg}{\mathsf{dg}}

\newcommand{\uHom}{\underline{\mathsf{Hom}}}
\newcommand{\Hom}{\mathsf{Hom}}

\newcommand{\HomC}{\uHom_{\,!}}
\newcommand{\HomA}{\uHom_{\mathsf{add}}}

\newcommand{\rep}{\mathsf{rep}}

\newcommand{\Cat}{\mathsf{Cat}}
\newcommand{\CAT}{\mathsf{CAT}}

\newcommand{\Ho}{\mathsf{Ho}}
\newcommand{\HO}{\mathsf{HO}}

\newcommand{\Hmo}{\mathsf{Hmo}}% Morita homotopy theory
\newcommand{\op}{\mathsf{op}}

\newcommand{\too}{\longrightarrow}

\newcommand{\HomL}{\uHom_{\mathsf{loc}}}

\newcommand{\ie}{\textsl{i.e.}\ }
\newcommand{\eg}{\textsl{e.g.}}

\newcommand{\Madd}{\Mot^{\mathsf{add}}_\dg}

\newcommand{\Uadd}{\cU^{\mathsf{add}}_{\dg}}

\newcommand{\Mloc}{\Mot^{\mathsf{loc}}_{\dg}}
\newcommand{\Uloc}{\cU^{\mathsf{loc}}_{\dg}}

\begin{document}

\title[A guided tour through the garden of noncommutative motives]{A guided tour through the \\garden of noncommutative motives}
\author{Gon{\c c}alo~Tabuada}

\address{Gon{\c c}alo Tabuada, Department of Mathematics, MIT, Cambridge, MA 02139}
\email{tabuada@math.mit.edu}

\date{\today}


%\thanks{The author was partially supported by the FCT-Portugal grants {\tt PTDC/MAT/098317/2008} and {\tt SFRH/BSAB/1116/2011}}

\abstract{These are the extended notes of a survey talk on noncommutative motives given at the {\it $3^{\mathrm{era}}$ Escuela de Inverno Luis Santal{\'o}-CIMPA Research School: Topics in Noncommutative Geometry}, Buenos Aires, July 26 to August 6, 2010.}
}

\maketitle
\vskip-\baselineskip
\vskip-\baselineskip
\vskip-\baselineskip
\tableofcontents
\vskip-\baselineskip
\vskip-\baselineskip

\noindent\textbf{Acknowledgments:} I would like to start by thanking all the organizers of this winter school and specially professor Guillermo Corti{\~n}as for kindly giving me the opportunity to present my work. I am also very grateful to the Clay Mathematics Institute for financially supporting my participation. A special thanks goes also to Eugenia Ellis who made me discover the magic of the {\em Argentine Tango} as well as the fascinating city of Buenos Aires. Finally, I would like to thank the Department of Mathematics of MIT for a warm reception.
$$\star \star \star $$
In order to make these notes accessible to a broad audience, I have decided to emphasize the conceptual ideas behind the theory of noncommutative motives rather than its technical aspects. I will start by stating two foundational questions. One concerning {\em higher algebraic $K$-theory} ({\bf Question A}) and another one concerning {\em noncommutative algebraic geometry} ({\bf Question B}). One of the main goals of this guided tour will be not only to provide precise answers to these distinct questions but moreover to explain what is the relation between the corresponding answers.

%-------------------------------------------------------------------------------
\section{Higher algebraic $K$-theory}
%-------------------------------------------------------------------------------
Algebraic K-theory goes back to Grothendieck's work \cite{SGA6} on the Riemann-Roch theorem. Given a commutative ring $R$ (or more generally an algebraic variety), he introduced the nowadays called {\em Grothendieck group} $K_0(R)$ of $R$. Latter, in the sixties, Bass~\cite{Bass} defined $K_1(R)$ as the abelianization of the general linear group $\mathrm{GL}(R)$. These two abelian groups, whose applications range from arithmetic to surgery of manifolds, are very well understood from a conceptual and computational point of view; see Weibel's survey \cite{Weibel-survey}. After Bass' work, it became clear that these groups should be part of a whole family of {\em higher} algebraic $K$-theory groups. After several attempts made by several mathematicians, it was Quillen who devised an elegant {\em topological} construction; see \cite{Quillen}. He introduced, the nowadays called Quillen's {\em plus construction} $(-)^+$, by which we simplify the fundamental group of a space without changing its (co-)homology groups. By applying this construction to the classifying space $\mathrm{BGL}(R)$ (where simplification in this case means abelianization), he defined the higher algebraic $K$-theory groups as
\begin{equation*}%\label{eq:mechanism}
K_n(R):= \pi_n(\mathrm{BGL}(R)^+ \times K_0(R)) \qquad n \geq 0\,.
\end{equation*}
Since Quillen's foundational work, higher algebraic $K$-theory has found extraordinary applications in a wide range of research fields; consult \cite{Handbook}. However, Quillen's mechanism for manufacturing these higher algebraic $K$-theory groups remained rather mysterious until today. Hence, the following question is of major importance:

\medbreak

{\bf Question A:} {\em How to conceptually characterize higher algebraic $K$-theory ?}

%-------------------------------------------------------------------------------
\section{Noncommutative algebraic geometry}
%-------------------------------------------------------------------------------
Noncommutative algebraic geometry goes back to Bondal-Kapranov's work \cite{BK,BK1} on exceptional collections of coherent sheaves. Since then, Drinfeld, Kaledin, Kontsevich, Orlov, Van den Bergh, and others, have made important advances; see \cite{Orlov,BB,Drinfeld,Chitalk,Kaledin,IAS,ENS, Miami,finMot}. Let $X$ be an algebraic variety. In order to study it, we can proceed in two distinct directions.

In one direction, we can associate to $X$ several (functorial) invariants like the Grothendieck group ($K_0$), the higher $K$-theory groups ($K_\ast$), the negative $K$-theory groups ($\bbK_\ast$), the cyclic homology groups ($HC_\ast$) and all its variants (Hochschild, periodic, negative, $\ldots$), the topological cyclic homology groups ($TC_\ast$), etc. Each one of these invariants encodes a particular arithmetic/geometric feature of the algebraic variety $X$.

In the other direction, we can associate to $X$ its derived category $\cD_\perf(X)$ of perfect complexes of $\cO_X$-modules. The importance of this triangulated category relies on the fact that any correspondance between $X$ and $X'$ which induces an equivalence between the derived categories $\cD_\perf(X)$ and $\cD_\perf(X')$, induces also an isomorphism on all the above invariants. Hence, it is natural to ask if the above invariants of $X$ can be recovered directly out of $\cD_\perf(X)$. This can be done in very particular cases (\eg\, the Grothendieck group) but not in full generality. The reason being is that when we pass from $X$ to $\cD_\perf(X)$ we loose too much information concerning $X$. We should therefore ``stop somewhere in the middle''. In order to formalize this insight, Bondal and Kapranov introduced the following notion.
\begin{definition}{(Bondal-Kapranov \cite{BK,BK1})}
A {\em differential graded (=dg) category} $\cA$, over a (fixed) base commutative ring $k$, is a category enriched over complexes of $k$-modules (morphism sets $\cA(x,y)$ are complexes) in such a way that composition fulfills the Leibniz rule: $d(f\circ g)=d(f)\circ g + (-1)^{\mathrm{deg}(f)}f \circ d(g)$. A {\em differential graded (=dg) functor} is a functor which preserves the differential graded structure; consult Keller's ICM adress~\cite{ICM} for further details. The category of (small) dg categories (over $k$) is denoted by $\dgcat$.
\end{definition}
Associated to the algebraic variety $X$ there is a natural dg category $\cD_\perf^\dg(X)$ which enhances\footnote{Consult Lunts-Orlov \cite{LO} for the uniqueness of this enhancement.} the derived category $\cD_\perf(X)$, \ie the latter category is obtained from the former one by applying the $0^{\mathrm{th}}$-cohomology group functor at each complex of morphisms. By considering $\cD_\perf^\dg(X)$ instead of $\cD_\perf(X)$ we solve many of the (technical) problems inherent to triangulated categories like the non-functoriality of the cone. More importantly, we are able to recover all the above invariants of $X$ directly out of $\cD_\perf^\dg(X)$. This circle of ideas is depicted in the following diagram:
$$
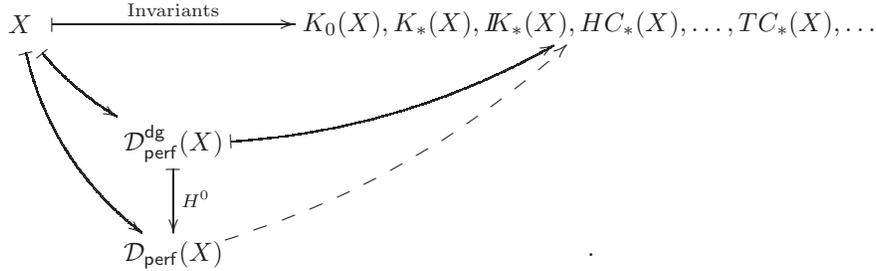
\xymatrix{
*+<3.5ex>{X} \ar@/_0.5pc/@{|->}[dr] \ar@/_1pc/@{|->}[ddr] \ar@{|->}[rr]^-{\mathrm{Invariants}} && K_0(X), K_\ast(X), \bbK_\ast(X), HC_\ast(X), \ldots, TC_\ast(X), \ldots \\
& \cD_\perf^\dg(X) \ar@{|->}[d]^{H^0} \ar@/_1pc/@{|->}[ur] & \\
& \cD_\perf(X) \ar@/_1pc/@{-->}[uur] & \,.\\
}
$$
From the point of view of the invariants, there is absolutely no difference between the algebraic variety $X$ and the dg category $\cD_\perf^\dg(X)$. This is the main idea behind noncommutative algebraic geometry: given a dg category, we should consider it as being the dg derived category of perfect complexes over a hypothetical noncommutative space and try to do ``algebraic geometry'' directly on it. Citing Drinfeld~\cite{Chitalk}, noncommutative algebraic geometry can be defined as:  {\it ``the study of dg categories and their homological invariants ''}.

\begin{example}{(Beilinson~\cite{Be})}\label{ex:Beilinson}
Suppose that $X$ is the $n^\mathrm{th}$-dimensional projective space $\bbP^n$. Then, there is an equivalence of dg categories
$$ \cD_\perf^\dg(\bbP^n) \simeq \cD_\perf^\dg(B)\,,$$
where $B$ is the algebra $\mathrm{End}(\cO(0) \oplus \cO(1) \oplus \ldots \oplus \cO(n))^\op$. Note that the abelian category of quasi-coherent sheaves on $\bbP^n$ is far from being the category of modules over an algebra. Beilinson's remarkable result show us that this situation changes radically when we pass to the derived setting. Intuitively speaking, the $n^{\mathrm{th}}$-dimensional projective space is an ``affine object'' in noncommutative algebraic geometry since it is described by a single (noncommutative) algebra.
\end{example}
In the commutative world, Grothendieck envisioned a theory of {\em motives} as a gateway between algebraic geometry and the assortment of the classical Weil cohomology theories (de Rham, Betti, $l$-adic, crystalline, and others); consult the monograph~\cite{Motives}.

In the noncommutative world we can envision a similiar picture. The role of the algebraic varieties and of the classical Weil cohomologies is played, respectively, by the dg categories and the numerous (functorial) invariants\footnote{In order to simplify the (graphical) exposition, we have decided to forget the $k$-linear structure of the cyclic homology groups $HC_\ast$.} 
\begin{equation}\label{eq:dia-main}
\xymatrix{
\dgcat \ar[rrrrrr]^-{K_\ast,\, \bbK_\ast,\, HC_\ast,\ldots, TC_\ast, \ldots} &&&&&& \mathrm{Ab}
}
\end{equation}
The Grothendieckian idea of motives consists then on combing this skein of invariants in order to isolate the truly fundamental one:
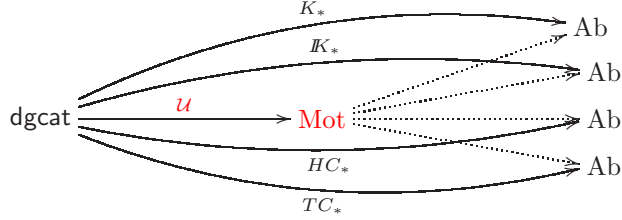
\begin{equation*}%\label{eq:dia-main1}
\vcenter{
\xymatrix@C=8em@R=0.4em{
&& \mathrm{Ab}
\\
&& \mathrm{Ab}  \kern-1em
\\
\dgcat \ar[r]^-{{\color{red}\cU}}  \ar@/^1.5pc/[rruu]^-{K_\ast}
 \ar@/^1pc/[rru]^-{\bbK_\ast} \ar@/_1pc/[rr]_-{HC_\ast}
 \ar@/_1.5pc/[rrd]_-{TC_\ast}& {\color{red} \Mot}
 \ar@{..>}[ruu]
 \ar@{..>}[ru] \ar@{..>}[r]
 \ar@{..>}[rd]
& \mathrm{Ab} \kern-1em
\\
&&\mathrm{Ab} \kern-1em
\\
}}\end{equation*}
The gateway category ${\color{red} \Mot}$, through which all invariants factor uniquely, should then be called the category of {\em noncommutative motives} and the functor ${\color{red} \cU}$ the {\em universal invariant}. Note that in this yoga, the different invariants are simply different representations of the motivic category ${\color{red} \Mot}$. In particular, any result which holds in ${\color{red} \Mot}$, holds everywhere. This beautiful circle of ideas leads us to the following down-to-earth question:

\medbreak

{\bf Question B:} {\em Is there a well-defined category of noncommutative motives ?}

%-------------------------------------------------------------------------------
\section{Derived Morita equivalences}
%-------------------------------------------------------------------------------
Note first that all the classical constructions which can be performed with $k$-algebras can also be performed with dg categories; consult \cite{ICM}. A dg functor $F:\cA \to \cB$ is called a {\em derived Morita equivalence} if the induced restriction of scalars functor $\cD(\cB) \stackrel{\sim}{\to} \cD(\cA)$ is an equivalence of (triangulated) categories. Thanks to the work of Blumberg-Mandell, Keller, Schlichting, and Thomason-Trobaugh, all the invariants \eqref{eq:dia-main} invert derived Morita equivalences; see \cite{BM,cyclic,Marco,TT}. Intuitively speaking, although defined at the ``dg level'', these invariants only depend on the underlying derived category. Hence, it is crucial to understand dg categories up to derived Morita equivalence. The following result is central in this direction. 
\begin{theorem}{(\cite{IMRN,cras})}\label{thm:Quillen}
The category $\dgcat$ carries a (cofibrantly generated) Quillen model structure\footnote{An analogous model structure in the setting of spectral categories was developed in \cite{Advances}.} whose weak equivalences are the derived Morita equivalences.
\end{theorem}
The homotopy category obtained is denoted by $\Hmo$. Theorem~\ref{thm:Quillen} allow us to study the purely algebraic setting of dg categories using ideas, techniques, and insights of topological nature. Here are some examples:
%-------------------------------------------------------------------------------
\subsection*{Bondal-Kapranov's pre-triangulated envelope}
%-------------------------------------------------------------------------------
Using ``one-sided twisted complexes'', Bondal and Kapranov constructed in \cite{BK} a pre-triangulated envelope $\cA^{\mathsf{pre}\text{-}\mathsf{tr}}$ of every dg category $\cA$. Intuitively speaking, their construction consists on formally adding to $\cA$ (de-)suspensions, cones, cones of morphisms between cones, etc. Thanks to Theorem~\ref{thm:Quillen}, this involved contribution can be conceptually characterized as being simply a functorial fibrant resolution functor; see \cite{IMRN}.
%-------------------------------------------------------------------------------
\subsection*{Drinfeld's DG quotient}
%-------------------------------------------------------------------------------
The most useful operation which can be performed on triangulated categories is the passage to a Verdier 
quotient. Recently, through a very elegant construction (reminiscent from Dwyer-Kan localization), Drinfeld \cite{Drinfeld} lifted this operation to the world of dg categories. Although very elegant, this construction didn't seem to satisfy any obvious universal property. Theorem~\ref{thm:Quillen} allowed us to complete this aspect of Drinfeld's work by characterizing the DG quotient as a homotopy cofiber construction; see \cite{DGquotient}.
%-------------------------------------------------------------------------------
\subsection*{Kontsevich's saturated dg categories}
%-------------------------------------------------------------------------------
Kontsevich understood precisely how to express smooth and properness in the noncommutative world.
\begin{definition}{(Kontsevich \cite{IAS,ENS})}\label{def:smooth}
A dg category $\cA$ is called:
\begin{itemize}
\item {\em smooth} if it is perfect as a bimodule over itself;
\item {\em proper} if its complexes of $k$-modules $\cA(x,y)$ are perfect;
\item {\em saturated} if it is smooth and proper.
\end{itemize}
\end{definition}
Definition~\ref{def:smooth} is justified by the following fact: given a quasi-compact and quasi-separated scheme $X$, the dg category $\cD_\perf^\dg(X)$ is smooth and proper if and only if $X$ is smooth and proper in the sense of classical algebraic geometry. Other examples of saturated dg categories appear in study of Deligne-Mumford stacks, quantum projective varieties, Landau-Ginzburg models, etc. 

Now, note that the tensor product of $k$-algebras extends naturally to dg categories. By deriving it (with respect to derived Morita equivalences), we obtain then a symmetric monoidal structure on $\Hmo$. Making use of it, the saturated dg categories can be conceptually characterized as being precisely the dualizable (or rigid) objects in the symmetric monoidal category $\Hmo$; see \cite{CT1}. As in any symmetric monoidal category, we can define the Euler characteristic of a dualizable object. In topology, for instance, the Euler characteristic of a finite $CW$-complex is the alternating sum of the number of cells. In $\Hmo$, we have the following result.  
\begin{proposition}{(Cisinski $\&$ Tab.~\cite{CT1})}\label{prop:Euler}
Let $\cA$ be a saturated dg category. Then, its Euler characteristic $\chi(\cA)$ in $\Hmo$ is the Hochschild homology\footnote{More generally, the trace of an endomorphisms is given by Hochschild homology with coefficients.} complex $HH(\cA)$ of $\cA$.
\end{proposition}
Proposition~\ref{prop:Euler} illustrates the Grothendieckian idea of combining the skein of invariants \eqref{eq:dia-main} ``as far as possible'' in order to understand, directly on ${\color{red} \Mot}$, their conceptual nature. By simply inverting the class of derived Morita equivalences, Hochschild homology can be conceptually understood as the Euler characteristic.

%-------------------------------------------------------------------------------
\section{Noncommutative pure motives}
%-------------------------------------------------------------------------------
In order to answer {\bf Question B} we need to start by identifying the properties common to all the invariants \eqref{eq:dia-main}. In the previous section we have already observed that they are {\em derived Morita invariant}, \ie they send derived Morita equivalences to isomorphisms. In this section, we identify another common property. An {\em upper triangular matrix} $M$ is given by 
$$
\begin{array}{rcl}
M & := & \begin{pmatrix} \cA & X \\ 0 &
  \cB \end{pmatrix}\,,
\end{array}
$$
where $\cA$ and $\cB$ are dg categories and $X$ is a $\cA\text{-}\cB$-bimodule. The totalization $|M|$ of $M$ is the dg category whose set of objects is the disjoint union of the sets of objects of $\cA$ and $\cB$, and whose morphisms are given by: $\cA(x,y)$ if $x,y \in \cA$; $\cB(x,y)$ if $x,y \in \cB$; $X(x,y)$ if $x \in \cA$ and $y \in \cB$; $0$ if $x \in \cB$ and $y \in \cA$. Composition is induced by the composition operation on $\cA$ and $\cB$, and by the $\cA\text{-}\cB$-bimodule structure of $X$. Note that we have two natural inclusion dg functors $\iota_\cA: \cA \to |M|$ and $\iota_\cB: \cB \to |M|$.
\begin{definition}\label{def:additive}
Let $E:\dgcat \to \mathsf{A}$ be a functor with values in an additive category. We say that $E$ is an {\em additive invariant of dg categories} if it is derived Morita invariant and satisfies the following condition: for every upper triangular matrix $M$, the inclusion dg functors $\iota_\cA$ and $\iota_\cB$ induce an isomorphism
$$ E(\cA)\oplus E(\cB) \stackrel{\sim}{\too} E(|M|)\,.$$
\end{definition}
It follows from the work of Blumberg-Mandell, Keller, Schlichting, and Thomason-Trobaugh, that all the invariants \eqref{eq:dia-main} satisfy additivity, and hence are additive invariant of dg categories; see \cite{BM,cyclic,Marco,TT}. The universal additive invariant of dg categories was constructed in \cite{IMRN}. It can be described\footnote{A similar construction in the setting of spectral categories was developed in \cite{IMRN-spectral}.} as follows: let $\Hmo_0$ be the category whose objects are the dg categories and whose morphisms are given by $\Hom_{\Hmo_0}(\cA,\cB):= K_0\rep(\cA,\cB)$, where $\rep(\cA,\cB)\subset \cD(\cA^\op \otimes^\bbL \cB)$ the full triangulated subcategory of those $\cA\text{-}\cB$-bimodules $X$ such that $X(a,-) \in \cD_\perf(\cB)$ for every object $a \in \cA$. Composition is induced by the tensor product of bimodules. Note that we have a natural functor
$$ \cU_{\mathsf{A}}: \dgcat \too \Hmo_0$$
which is the identity on objects and which maps a dg functor to the class (in the Grothendieck group) of the naturally associated bimodule. The category $\Hmo_0$ is additive and the functor $\cU_{\mathsf{A}}$ is additive in the sense of Definition~\ref{def:additive}. Moreover, it is characterized by the following universal property.
\begin{theorem}{(\cite{IMRN})}\label{thm:additive}
Given an additive category $\mathsf{A}$, we have an induced equivalence of categories
\begin{equation*}%\label{eq:equiv1}
(\cU_\mathsf{A})^\ast: \mathrm{Fun}_{\mathsf{add}}(\Hmo_0,\mathsf{A}) \stackrel{\sim}{\too} \mathrm{Fun}_{\mathsf{additivity}}(\dgcat, \mathsf{A})\,,
\end{equation*}
where the left hand-side denotes the category of additive functors and the right hand-side the category of additive invariants in the sense of Definition~\ref{def:additive}.
\end{theorem}
The additive category $\Hmo_0$ (and $\cU_{\mathsf{A}}$) is our first answer to {\bf Question B}. A second answer will be described in Section~\ref{sec:NCmixed}. Note that by Theorem~\ref{thm:additive}, all the invariants \eqref{eq:dia-main} factor uniquely through $\Hmo_0$. This motivic category enabled several (tangential) applications. Here are two examples:
\begin{example}{(Chern characters)} The Chern character maps are one of the most important working tools in mathematics. Although they admit numerous different constructions, they were not fully understood at the conceptual level. Making use of the additive category $\Hmo_0$ and of Theorem~\ref{thm:additive} we have bridged this gap by characterizing the Chern character maps, from the Grothendieck group to the (negative) cyclic homology groups, in terms of simple universal properties; see \cite{Chern}.
\end{example}
\begin{example}{(Fundamental theorem)} The fundamental theorems in homotopy algebraic $K$-theory and periodic cyclic homology, proved respectively by Weibel~\cite{Weibel-KH} and Kassel~\cite{Kassel}, are of major importance. Their proofs are not only very different but also quite involved. Making use of the additive category $\Hmo_0$ and of Theorem~\ref{thm:additive}, we have given a simple, unified and conceptual proof of these fundamental theorems; see \cite{Fund}. 
\end{example}
%-------------------------------------------------------------------------------
\subsection*{Noncommutative Chow motives}
%-------------------------------------------------------------------------------
By restricting himself to saturated dg categories, which morally are the ``noncommutative smooth projective varieties'', Kontsevich introduced the following category.
\begin{definition}{(Kontsevich \cite{IAS,finMot}; \cite{CvsNC})}
Let $F$ be a field of coefficients. The category $\NChow_F$ of {\em noncommutative Chow motives} (over the base ring $k$ and with coefficients in $F$) is defined as follows: first consider the $F$-linearization $\Hmo_{0;F}$ of the additive category $\Hmo_0$. Then, pass to its idempotent completion $\Hmo_{0;F}^\natural$. Finally, take the idempotent complete full subcategory of $\Hmo_{0;F}^\natural$ generated by the saturated dg categories.
\end{definition}
The precise relation between the classical category of Chow motives and the category of noncommutative Chow motives is the following: recall that the category $\Chow_\bbQ$ of $\Chow$ motives (with rational coefficients) is $\bbQ$-linear, additive and symmetric monoidal. Moreover, it is endowed with an important $\otimes$-invertible object, namely the Tate motive $\bbQ(1)$. The functor $-\otimes \bbQ(1)$ is an automorphism of $\Chow_\bbQ$ and so we can consider the associated orbit category $\Chow(k)_\bbQ/_{\!\!-\otimes \bbQ(1)}$; consult \cite{CvsNC} for details. Informally speaking, Chow motives which differ from a Tate twist become isomorphic in the orbit category.
\begin{theorem}{(Kontsevich \cite{IAS,finMot}; \cite{CvsNC}})\label{thm:Chow}
There exists a fully-faithful, $\bbQ$-linear, additive, and symmetric monoidal functor $R$ making the diagram 
\begin{equation}\label{eq:diagram}
\xymatrix@C=2em@R=1.5em{
\mathsf{SmProj}^\op \ar[rr]^{\cD_\perf^\dg(-)} \ar[d]_M && \dgcat \ar[d]^{\cU_{\mathsf{A}}} \\
\Chow_\bbQ \ar[d]_\pi && \Hmo_0 \ar[d]^{(-)_\bbQ^\natural} \\
\Chow_\bbQ\!/_{\!\!-\otimes \bbQ(1)} \ar[rr]_-R && \NChow_\bbQ \subset \Hmo_{0;\bbQ}^\natural \qquad \qquad
}
\end{equation}
commute (up to a natural isomorphism).
\end{theorem}
Intuitively speaking, Theorem~\ref{thm:Chow} formalizes the conceptual idea that the commutative world can be embedded into the noncommutative world after factorizing out by the action of the Tate motive. The above diagram \eqref{eq:diagram} opens new horizonts and opportunities of research by enabling the interchange of results, techniques, ideas, and insights between the commutative and the noncommutative world. This yoga was developed in \cite{CvsNC} in what regards Schur and Kimura finiteness, motivic measures, and motivic zeta functions. 

%-------------------------------------------------------------------------------
\subsection*{Noncommutative numerical motives}
%-------------------------------------------------------------------------------
In order to formalize and solve ``counting problems'', such as counting the number of common points to two planar curves in general position, the classical category of Chow motives is not appropriate as it makes use of a very refined notion of equivalence. Motivated by these ``counting problems'', Grothendieck developed in the sixties the category $\Num_F$ of numerical motives; see \cite{Motives}. Its noncommutative analogue can be described as follows: let $\cA$ and $\cB$ be two saturated dg categories and $\underline{X}=[\sum_i a_i X_i]\in \Hom_{\NChow_F}(\cA,\cB)$ and $\underline{Y}=[\sum_j b_j Y_j] \in \Hom_{\NChow_F}(\cB,\cA)$ two {\em noncommutative correspondances}. Their {\em intersection number} is given by the formula
\begin{equation}\label{eq:intersection1}
\langle \underline{X} \cdot \underline{Y}\rangle := \sum_{i,j,n} (-1)^n\, a_i\!\cdot \!b_j \!\cdot \!\mathrm{rk} HH_n(\cA; X_i \otimes^\bbL_\cB Y_j) \in F\,,
\end{equation}
where $\mathrm{rk} HH_n(\cA; X_i \otimes_\cB Y_j)$ denotes the rank of the $n^\mathrm{th}$-dimensional Hochschild homology group of $\cA$ with coefficients in the $\cA\text{-}\cA$-bimodule $X_i \otimes^\bbL_\cB Y_j$. A noncommutative correspondance $\underline{X}$ is {\em numerically equivalent to zero} if for every noncommutative correspondence $\underline{Y}$ the intersection number $\langle \underline{X} \cdot \underline{Y} \rangle$ is zero. As proved in \cite{Semi}, these correspondences form a $\otimes$-ideal of $\NChow(k)_F$, which we denote by $\cN$. 
\begin{definition}{({\it Marcolli $\&$ Tab.}~\cite{Semi})}\label{def:NCNum}
The category $\NNum_F$ of {\em noncommutative numerical motives} (over the base ring $k$ and with coefficients in $F$) is the idempotent completion of the quotient category $\NChow_F/\cN$. 
\end{definition}
The relation between Chow motives and noncommutative motives described in diagram~\eqref{eq:diagram} admits the following numerical analogue.
\begin{theorem}{(Marcolli $\&$ Tab.~\cite{Semi})}\label{thm:embedding}
There exists a fully-faithful, $\bbQ$-linear, additive, and symmetric monoidal functor $R_{\cN}$ making the diagram
$$
\xymatrix@C=2em@R=1.5em{
& \Chow_\bbQ  \ar[d]^\pi \ar[dl] &&  \\
\Num_\bbQ \ar[d]_\pi & \Chow_\bbQ/_{\!\!-\otimes \bbQ(1)}  \ar[dl]  \ar[rr]^-R && \ar[dl] \NChow_\bbQ  \\
\Num_\bbQ/_{\!\!-\otimes \bbQ(1)} \ar[rr]_{R_{\cN}} & & \NNum_\bbQ & 
}
$$
commute (up to natural isomorphism).
\end{theorem}
Intuitively speaking, Theorem~\ref{thm:embedding} formalizes the conceptual idea that Hochschild homology is the correct way to express ``counting'' in the noncommutative world. In the commutative world, Grothendieck conjectured that the category of numerical motives $\Num_F$ was abelian semi-simple. Jannsen~\cite{Jannsen}, thirty years latter, proved this conjecture without the use of any of the standard conjectures. Recently, we gave a further step forward by proving that Grothendieck's conjecture holds more broadly in the noncommutative world.
\begin{theorem}{(Marcolli $\&$ Tab.~\cite{Semi})}\label{thm:semi-simple}
Assume one of the following two conditions:
\begin{itemize}
\item[(i)] The base ring $k$  is local (or more generally that $K_0(k)=\bbZ$) and $F$ is a $k$-algebra; a large class of examples is given by taking $k=\bbZ$ and $F$ an arbitrary field.
\item[(ii)] The base ring $k$ is a field extension of $F$; a large class of examples is given by taking $F=\bbQ$ and $k$ a field of characteristic zero.
\end{itemize}
Then, the category $\NNum_F$ is abelian semi-simple. Moreover, if $\cJ$ is a $\otimes$-ideal in $\NChow_F$ for which the idempotent completion of the quotient category $\NChow_F/\cJ$ is abelian semi-simple, then $\cJ$ agrees with $\cN$.
\end{theorem}
Roughly speaking, Theorem~\ref{thm:semi-simple} shows that the unique way to obtain an abelian semi-simple category out of $\NChow_F$ is through the use of the above ``counting formula'' \eqref{eq:intersection1}, defined in terms of Hochschild homology. Among other applications, Theorem~\ref{thm:semi-simple} allowed us to obtain an alternative proof of Jannsen's result; see \cite{Semi}.
%-------------------------------------------------------------------------------
\subsection*{Kontsevich's noncommutative numerical motives}
%-------------------------------------------------------------------------------
Making use of a well-behaved bilinear form on the Grothendieck of saturated dg categories, Kontsevich introduced in \cite{IAS} a category $\mathrm{NCNum}_F$ of noncommutative numerical motives. Via duality arguments, the authors proved the following agreement result.
\begin{theorem}{(Marcolli $\&$ Tab.~\cite{Konts})}\label{thm:equality}  
The categories $\mathrm{NCNum}_F$ and $\NNum_F$ are equivalent.
\end{theorem}
By combining Theorem~\ref{thm:equality} with Theorem~\ref{thm:semi-simple}, we then conclude that $\mathrm{NCNum}_F$ is abelian semi-simple. Kontsevich conjectured this latter result in the particular case where $F=\bbQ$ and $k$ is of characteristic zero. We observe that Kontsevich's beautiful insight not only holds much more generally, but moreover it does not require the assumption of any (polarization) conjecture.
%-------------------------------------------------------------------------------
\section{Noncommutative mixed motives}\label{sec:NCmixed}
%-------------------------------------------------------------------------------
Up to now, we have been considering invariants with values in additive categories. From now on we will consider ``richer invariants'', taking values not in additive categories but in ``highly structured'' triangulated categories. In order to make this precise we will use the language of {\em Grothendieck derivators}, a formalism which allow us to state and prove precise universal properties; the reader who is unfamiliar with this language is invited to consult Appendix~\ref{appendix:A} at this point. Recall from Drinfeld~\cite{Drinfeld} that a sequence of dg functors $\cA \stackrel{I}{\to} \cB \stackrel{P}{\to} \cC$ is called {\em exact} if the induced sequence of derived categories $\cD(\cA)\to \cD(\cB) \to \cD(\cC)$ is exact in the sense of Verdier. For example, if $X$ is quasi-compact and quasi-separated scheme, $U \subset X$ a quasi-compact open subscheme and $Z:=X\backslash U$ the closed complementary, then the sequence of dg functors 
$$  \cD_\perf^\dg(X)_Z \too \cD_\perf^\dg(X) \too \cD_\perf^\dg(U)$$
is exact; see~Thomason-Trobaugh \cite{TT}. An exact sequence of dg functors is called {\em split-exact} if there exist dg functors $R:\cB \to \cA$ and $S: \cC \to \cB$, right adjoints to $I$ and $P$, respectively, such that $R \circ I \simeq \Id$ and $P \circ S \simeq \Id$ via the adjunction morphisms; consult \cite{Duke} for details.
\begin{definition}
Let $E: \HO(\dgcat) \to \bbD$ be a filtered homotopy colimit preserving morphism of derivators, from the derivator associated to the Quillen model structure of Theorem~\ref{thm:Quillen}, towards a strong triangulated derivator. We say that $E$ is a {\em localizing invariant} if it sends exact sequences to distinguished triangles 
\begin{eqnarray*}
\cA \too \cB \too \cC & \mapsto & E(\cA) \too E(\cB) \too E(\cC) \too E(\cA)[1]
\end{eqnarray*}
in the base category $\bbD(e)$ of $\bbD$. We say that $E$ is an {\em additive invariant} if it sends split-exact sequences to direct sums
\begin{eqnarray*}
\xymatrix@C=1.5em@R=1.0em{
  \cA \ar[r]  & \cB \ar[r]  \ar@/_0.5pc/[l] & \cC \ar@/_0.5pc/[l] 
} &\mapsto & E(\cA) \oplus E(\cC) \stackrel{\sim}{\to} E(\cB)\,.
\end{eqnarray*}
\end{definition}
Clearly, every localizing invariant is additive. Here are some classical examples.
\begin{example}{(Connective $K$-theory)}\label{ex:connec}
As explained in \cite{Duke}, connective $K$-theory gives rise to an additive invariant
$$K:\HO(\dgcat) \too \HO(\Spt)$$
with values in the triangulated derivator associated to the (stable) Quillen model category of spectra. Quillen's higher $K$-theory groups $K_\ast$ can then be obtained from this spectrum by taking stable homotopy groups. This invariant, although additive, is {\em not} localizing. The following example corrects this default. 
\end{example}
\begin{example}{(Nonconnective $K$-theory)}
As explained in \cite{Duke}, nonconnective $K$-theory gives rise to a localizing invariant
$$\bbK:\HO(\dgcat) \too \HO(\Spt)\,.$$
As in the previous example, Bass' negative algebraic $K$-theory groups $\bbK_\ast$ can be obtained from this spectrum by taking (negative) stable homotopy groups.
\end{example}
\begin{example}{(Mixed complex)}
Following Kassel \cite{Kassel}, let $\Lambda$ be the dg algebra $k[\epsilon]/\epsilon^2$ where $\epsilon$ is of degree $-1$ and $d(\epsilon)=0$. Under this notation, a {\em mixed complex} is simply a right dg $\Lambda$-module. As explained in \cite{Duke}, the mixed complex construction gives rise to a localizing invariant
$$ C: \HO(\dgcat) \too \HO(\Lambda\text{-}\mathrm{Mod})$$
with values in the triangulated derivator associated to the (stable) Quillen model category of right dg $\Lambda$-modules. Cyclic homology and all its variants (Hochschild, periodic, negative, $\ldots$) can be obtained from this mixed complex construction by a simple procedures; see \cite{Kassel}.
\end{example}
\begin{example}{(Topological cyclic homology)}\label{ex:TC}
As explained by Blumberg and Mandell in \cite{BM} (see also \cite{AGT}), topological cyclic homology gives rise to a localizing invariant
$$ TC: \HO(\dgcat) \too \HO(\Spt)\,.$$
The topological cyclic homology groups $TC_\ast$ can be obtained from this spectrum by taking stable homotopy groups.
\end{example}
In order to simultaneously study all the above classical examples, the universal additive and localizing invariants 
\begin{eqnarray*}%\label{eq:morphisms}
\Uadd: \HO(\dgcat) \too \Madd && \Uloc: \HO(\dgcat) \too \Mloc
\end{eqnarray*}
were constructed\footnote{A similar approach in the setting of $\infty$-categories was developed by Blumberg, Gepner and the author in \cite{infinity}. Besides algebraic and geometric examples, the authors studied also topological examples like $A$-theory.} in \cite{Duke}. They are characterized (in the $2$-category of Grothendieck derivators) by the following universal property.
\begin{theorem}{(\cite{Duke})}\label{thm:main2}
Given a strong triangulated derivator $\bbD$, we have induced equivalences of categories
\begin{eqnarray*}
(\Uadd)^\ast : \uHom_!(\Madd, \bbD) & \stackrel{\sim}{\too} &  \HomA(\HO(\dgcat), \bbD)\\
(\Uloc)^\ast : \uHom_!(\Mloc, \bbD) & \stackrel{\sim}{\too} & \HomL(\HO(\dgcat), \bbD)\,,
\end{eqnarray*}
where the right hand-sides denote, respectively, the categories of additive and localizing invariants.
\end{theorem}
\begin{remark}{(Quillen model)}
The additive and the localizing motivator admit natural Quillen models given in terms of a Bousfield localization of presheaves of (symmetric) spectra; consult \cite{Duke} for details.
\end{remark}
Because of these universal properties, $\Madd$ is called the {\em additive motivator}, $\Mloc$ the {\em localizing motivator}, $\Uadd$ the {\em universal additive invariant}, $\Uloc$ the {\em universal localizing invariant}, $\Madd(e)$ the {\em triangulated category of noncommutative additive motives}, and $\Mloc(e)$ the {\em triangulated category of noncommutative localizing motives}. Note that since localization implies additivity, we have a well-defined (homotopy colimit preserving) morphism of derivators $\Madd \to \Mloc$. The triangulated category $\Madd(e)$ (and $\Mloc(e)$) is our second answer to {\bf Question B}. Note that by Theorem~\ref{thm:main2}, all the invariants of Examples \ref{ex:connec}-\ref{ex:TC} factor uniquely through $\Madd(e)$. Since the composed functor
$$ \dgcat \too \Hmo \stackrel{\Uadd(e)}{\too} \Madd(e)$$
is an additive invariant of dg categories in the sense of Definition~\ref{def:additive}, we obtain by Theorem~\ref{thm:additive} an induced additive functor $\Hmo_0 \to \Madd(e)$, which turns out to be fully-faithful. Intuitively speaking, our second answer to {\bf Question B} contains the first one. In other words, the world of noncommutative pure motives is contained in the world of noncommutative mixed motives. As we will see in the next section, the latter world is much richer than the former one.

In Example~\ref{ex:Beilinson}, we observed that the dg category $\cD_\perf^\dg(\bbP^n)$ is derived Morita equivalent to the algebra $\mathrm{End}(\cO(0)\oplus \cO(1) \oplus \ldots \oplus \cO(n))^\op$. By passing to the triangulated category of noncommutative additive motives, we obtain the following splitting:
$$ \Uadd(\cD_\perf^\dg(\bbP^n))\simeq \underbrace{\Uadd(k) \oplus \cdots \oplus \Uadd(k)}_{(n+1)\text{-}\mathrm{copies}}\,.$$ 
The reason behind this phenomenon is a semi-orthogonal decomposition of the triangulated category $\cD_\perf(X)$. Intuitively speaking, the noncommutative additive motive of the $n^{\mathrm{th}}$-dimensional projective space consists simply of $n+1$ ``points''. 

The motivic category $\Madd(e)$ enabled several (tangential) applications, Here is one illustrative example:
\begin{example}{(Farrell-Jones isomorphism conjectures)}
The Farrell-Jones isomorphism conjectures are important driving forces in current mathematical research and imply well-know conjectures due to Bass, Borel, Kaplansky, Novikov; see L{\"u}ck-Reich's survey in \cite{Handbook}. Given a group $G$, they predict the value of algebraic $K$- and $L$-theory of the group ring $k[G]$ in terms of its values on the virtually cyclic subgroups of $G$. In addition, the literature contains many variations of this theme, obtained by replacing the $K$- and $L$-theory functors by other functors like Hochschild homology, topological
cyclic homology, etc. During the last decades each one of these isomorphism conjectures has been
proved for large classes of groups using a variety of different methods.
Making use of Theorem~\ref{thm:main2}, Balmer and the author organized this exuberant herd of conjectures by
explicitly describing the fundamental isomorphism conjecture; see \cite{BT}. It turns out that this fundamental
conjecture, which implies all the existing isomorphism conjectures on the market, can be described
solely in terms of algebraic $K$-theory. More precisely, it is a simple ``coefficient variant'' of the classical
Farrell-Jones conjecture in algebraic $K$-theory.
\end{example}
%-------------------------------------------------------------------------------
\section{Co-representability}
%-------------------------------------------------------------------------------
As in any triangulated derivator, the additive and localizing motivators are canonically enriched over spectra. Let us denote by $\bbR\Hom(-,-)$ their spectra of morphisms; see Appendix~\ref{appendix:A}. Connective algebraic $K$-theory is an example of an additive invariant while nonconnective algebraic $K$-theory is an example of a localizing invariant. Therefore, by Theorem~\ref{thm:main2}, they descend to the additive and localizing motivator, respectively. The following result show us that they become co-representable by the noncommutative motive associated to the base ring.
\begin{theorem}{(\cite{Duke}; Cisinski $\&$ Tab.~\cite{CT})}\label{thm:2}
Given a dg category $\cA$, we have natural equivalences of spectra
\begin{eqnarray*}
\bbR\Hom(\Uadd(k),\Uadd(\cA)) \simeq K(\cA) && \bbR\Hom(\Uloc(k),\Uloc(\cA)) \simeq \bbK(\cA) \,.
\end{eqnarray*}
In the triangulated categories of noncommutative motives, we have natural isomorphisms of abelian groups
\begin{eqnarray*}
\Hom(\Uadd(k), \Uadd(\cA)[-n]) & \simeq & K_n (\cA) \qquad n \geq 0 \\
\Hom(\Uloc(k), \Uloc(\cA)[-n]) & \simeq & \bbK_n( \cA) \qquad n \in \bbZ\,.
\end{eqnarray*}
\end{theorem}
\begin{example}{(Schemes)}
By taking $\cA=\cD_\perf^\dg(X)$ in Theorem~\ref{thm:2}, with $X$ a quasi-compact and quasi-separated scheme, we recover the connective $K(X)$ and nonconnective $\bbK(X)$ $K$-theory spectrum of $X$.
\end{example}
\begin{remark}{(Bivariant $K$-theory)}
Theorem~\ref{thm:2} is in fact richer. In what concerns the additive motivator, the base ring $k$ can be replaced by any homotopically finitely presented dg category $\cB$ (the homotopical version of the classical notion of finite presentation) and $K(\cA)$ by the bivariant $K$-theory of $\cB\text{-}\cA$-bimodules. In what concerns the localizing motivator, the base ring $k$ can be replaced by any saturated dg category $\cB$ and $\bbK(\cA)$ by the spectrum $\bbK(\cB^\op \otimes \cA)$; consult \cite{CT,CT1,Duke}.
\end{remark}

\begin{remark}{(Bivariant cyclic homology)}
Classical theories like bivariant cyclic cohomology (and the associated Connes' bilinear pairings) can also be expressed as morphisms sets in the category of noncommutative motives; see \cite{Bivariant}.
\end{remark}

Theorem~\ref{thm:2} is our answer to {\bf Question A}. Note that while the right-hand sides are, respectively, connective and nonconnective algebraic $K$-theory, the left-hand sides are defined solely in terms of precise universal properties: algebraic $K$-theory is never used (or even mentioned) in their construction. Hence, the equivalences of Theorem~\ref{thm:2} provide us with a conceptual characterization of higher algebraic $K$-theory. To the best of the author's knowledge, this is the first conceptual characterization of algebraic $K$-theory since Quillen's foundational work. We can even take these equivalences as the very definition of higher algebraic $K$-theory. The precise relation between the answers to {\bf Questions A} and {\bf B} is by now clear. Intuitively speaking, connective (resp. nonconnective) algebraic $K$-theory is the additive (resp. localizing) invariant co-represented by the noncommutative motive associated to the base ring, which as explained in the next section is simply the $\otimes$-unit object.

%-------------------------------------------------------------------------------
\section{Symmetric monoidal structure}
%-------------------------------------------------------------------------------
The tensor product of $k$-algebras extends naturally to dg categories, giving rise to a symmetric monoidal structure on $\HO(\dgcat)$. The $\otimes$-unit is the base ring $k$ (considered as a dg category). Making use of a derived version of Day's convolution product, the authors proved the following result.
\begin{theorem}{(Cisinski $\&$ Tab.~\cite{CT1})}\label{thm:sym}
The additive and localizing motivators carry a canonical symmetric monoidal structure making the universal additive and localizing invariants symmetric monoidal. Moreover, these symmetric monoidal structures preserve homotopy colimits in each variable and are characterized by the following universal property: given any strong triangulated derivator $\bbD$, endowed with a symmetric monoidal structure, we have induced equivalence of categories:
\begin{eqnarray*}
(\Uadd)^\ast : \uHom^\otimes_!(\Madd, \bbD) & \stackrel{\sim}{\too} &  \HomA^\otimes(\HO(\dgcat), \bbD)\\
(\Uloc)^\ast : \uHom^\otimes_!(\Mloc, \bbD) & \stackrel{\sim}{\too} & \HomL^\otimes(\HO(\dgcat), \bbD)\,.
\end{eqnarray*}
\end{theorem}
%-------------------------------------------------------------------------------
\subsection*{Kontsevich's noncommutative mixed motives}
%-------------------------------------------------------------------------------
In \cite{IAS,finMot}, Kontsevich introduced a category $\Mix$ of {\em noncommutative mixed motives} (over the base ring $k$). Roughly speaking, $\Mix$ is obtained by taking a formal idempotent completion of the triangulated envelope of the category of saturated dg categories (with bivariant algebraic $K$-theory spectra as morphism sets). Making use Theorem~\ref{thm:sym}, the category $\Mix$ can be ``realized'' inside the triangulated category of noncommutative motives.
\begin{proposition}{(Cisinski $\&$ Tab.~\cite{CT1})}\label{prop:embedding}
There is a natural fully-faithful embedding (enriched over spectra) of Kontsevich's category $\Mix$ of noncommutative mixed motives into the triangulated category $\Mloc(e)$ of noncommutative localizing motives. The essential image is the thick triangulated subcategory spanned by the noncommutative motives of saturated dg categories.
\end{proposition}
%Note that Proposition~\ref{prop:embedding} combined with the universality of $\Mloc(e)$ justifies the correctness of Kontsevich's {\em ad hoc} construction.
\begin{remark}{(Relation with Voevodsky's motives)}
In the same vein as Theorem~\ref{thm:Chow}, Voevodsky's triangulated category $\mathrm{DM}$ of motives~\cite{Voevodsky} relates to (a $\bbA^1$-homotopy variant of) Kontsevich's category $\Mix$ of noncommutative mixed motives. The author and Cisinski are nowadays in the process of writing down this result.
\end{remark}
%-------------------------------------------------------------------------------
\subsection*{Products in algebraic $K$-theory}
%-------------------------------------------------------------------------------
Let $\cA$ and $\cB$ be two dg categories. On one hand, following Waldhausen \cite{Wald}, we have a classical algebraic $K$-theory pairing
\begin{eqnarray}\label{eq:pairing1}
K(\cA) \wedge K(\cB) \too K(\cA \otimes \cB)\,.
\end{eqnarray}
On the other hand, by combining the co-representability Theorem~\ref{thm:2} with Theorem~\ref{thm:sym}, we obtain another well-defined algebraic $K$-theory pairing
\begin{eqnarray}\label{eq:pairing2}
K(\cA) \wedge K(\cB) \too K(\cA \otimes \cB)\,.
\end{eqnarray}
\begin{theorem}{(\cite{Prods})}\label{thm:prods}
The pairings \eqref{eq:pairing1} and \eqref{eq:pairing2} agree up to homotopy; a similar result holds for nonconnective $K$-theory.
\end{theorem}
\begin{example}{(Commutative algebras)}\label{ex:prods}
Let $\cA=\cB=A$, with $A$ is a {\em commutative} $k$-algebra. Then, by composing the pairing \eqref{eq:pairing2} with the multiplication map 
$$K(A\otimes A) \simeq \bbR\Hom(\Uadd(k), \Uadd(A \otimes A)) \too \bbR\Hom(\Uadd(k), \Uadd(A)) \simeq K(A)$$ 
we recover inside $\Madd$ the algebraic $K$-theory pairing on $K(A)$ constructed originally by Waldhausen. In particular, we recover the (graded-commutative) multiplicative structure on $K_\ast(A)$ constructed originally by Loday~\cite{Loday}. 
\end{example}
\begin{example}{(Schemes)}\label{ex:prods1}
When $\cA=\cB=\cD_\perf^\dg(X)$, with $X$ a quasi-compact and quasi-separated $k$-scheme, an argument similar to the one of the above example allow us to recover inside $\Madd$ the algebraic $K$-theory pairing on $X$ constructed originally by Thomason-Trobaugh~\cite{TT}.
\end{example}
Theorem~\ref{thm:prods} (and Examples \ref{ex:prods}-\ref{ex:prods1}) offers an elegant conceptual characterization of the algebraic $K$-theory products. Intuitively speaking, while Theorem~\ref{thm:2} shows us that connective algebraic $K$-theory is the additive invariant co-represented by the $\otimes$-unit of $\Madd$, Theorem~\ref{thm:prods} shows us that the classical algebraic $K$-theory products are simply the operations naturally induced by the symmetric monoidal structure on $\Madd$. 

%-------------------------------------------------------------------------------
\section{Higher Chern characters}
%-------------------------------------------------------------------------------
Higher algebraic $K$-theory is a very powerful and subtle invariant whose calculation is often out of reach. In order to capture some of its information, Connes-Karoubi, Dennis, Goodwillie, Hood-Jones, Kassel, McCarthy, and others, constructed higher Chern characters towards simpler theories by making use of a variety of highly involved techniques; see \cite{CK,Dennis,Goodwillie,HJ,Kassel-biv,McCarthy}. 

Making use of the theory of noncommutative motives, these higher Chern characters can be constructed, and conceptually characterized, in a simple and elegant way; see \cite{CT,CT1,Duke,Prods,AGT}. Let us now illustrate this in a particular case: choose your favorite additive invariant $E$ with values in the derivator associated to spectra. A classical example is given by connective algebraic $K$-theory. Thanks to Theorem~\ref{thm:main2}, we obtain then (homotopy colimit preserving) morphisms of derivators
$$ \overline{K}, \overline{E}: \Madd \too \HO(\Spt)$$
such that $\overline{K} \circ \Uadd = K$ and $\overline{E} \circ \Uadd = E$. Recall from Theorem~\ref{thm:2} that the functor $\overline{K}$ is co-represented by the noncommutative additive motive $\Uadd(k)$. Hence, the enriched Yoneda lemma furnishes us a natural equivalence of spectra $\bbR\Nat(\overline{K},\overline{E}) \simeq \overline{E}(k)$, where $\bbR\Nat$ denotes the spectrum of natural transformations. Using Theorem~\ref{thm:main2} again, we obtain a natural equivalence $\bbR\Nat(K,E) \simeq E(k)$. By passing to the $0^{\mathrm{th}}$-homotopy group, we conclude that there is a natural bijection between the natural transformation (up to homotopy) from $K$ to $E$ and $\pi_0E(k)$. In sum, the theory of noncommutative motives allow us to fully classify in simple and elegant terms all possible natural transformation from connective $K$-theory towards any additive invariant; a similar result holds for nonconnective $K$-theory.
\begin{example}{(Chern character)}\label{ex:HC}
Let $E$ be the cyclic homology $HC$ additive functor (promoted to an invariant taking values in spectra). Then, we have the following identifications:
\begin{eqnarray*}
\Nat(K,HC) \stackrel{\sim}{\to} k\simeq HC_0(k) && \{\mathrm{Chern}\,\,\mathrm{character}\} \mapsto {\bf 1}\,.
\end{eqnarray*}
\end{example}
Example~\ref{ex:HC} provides a conceptual characterization of the Chern character as being precisely the unit among all possible natural transformations. A similar characterization of the cyclotomic trace map, in the setting of $\infty$-categories, was recently developed by Blumberg, Gepner and the author in \cite{trace}. 

\appendix

%-----------------------------------------------------------------------
\section{Grothendieck derivators }\label{appendix:A}
%-----------------------------------------------------------------------
The original reference for the theory of derivators is Grothendieck's original manuscript~\cite{Grothendieck}. See also a short account by Cisinski and Neeman in \cite{CN}. Derivators originate in the problem of higher homotopies in derived categories.
For a triangulated category $\cT$ and for $X$ a
small category, it essentially never happens that the diagram
category $\Fun(X,\cT)=\cT^X$ remains triangulated; it already fails for
the category of arrows in~$\cT$, that is, for
$X=(\bullet\to\bullet)$. Now, very often, our triangulated category $\cT$ appears as the
homotopy category $\cT=\Ho(\cM)$ of some Quillen model~$\cM$. In this case, we can consider the category
$\Fun(X,\cM)$ of diagrams in $\cM$, whose homotopy category
$\Ho(\Fun(X,\cM))$ is often triangulated and provides a reasonable
approximation for $\Fun(X,\cT)$. More importantly, one can let
$X$ move. This nebula of categories $\Ho(\Fun(X,\cM))$, indexed by
small categories~$X$, and the various functors and natural
transformations between them is what Grothendieck formalized into
the concept of \emph{derivator}.

A derivator $\bbD$ consists of a
strict contravariant $2$-functor from the $2$-category of small
categories to the $2$-category of all categories
$$
\bbD: \Cat^\op \longrightarrow \CAT,
$$
subject to certain conditions; consult~\cite{CN} for details. The
essential example to keep in mind is the derivator $\bbD=\HO(\cM)$
associated to a (cofibrantly generated) Quillen model category~$\cM$
and defined for every small category~$X$ by
\begin{equation*}%\label{eq:HO}
\HO(\cM)(X)=\Ho(\Fun(X^\op,\cM))\,.
\end{equation*}
We denote by $e$ the $1$-point category with one object and one
identity morphism. Heuristically, the category $\bbD(e)$ is the
basic ``derived" category under consideration in the
derivator~$\bbD$. For instance, if $\bbD=\HO(\cM)$ then
$\bbD(e)=\Ho(\cM)$. Let us now recall two slightly technical properties of derivators.
\begin{enumerate}
\item[-]
A derivator $\bbD$ is called {\em strong} if for every finite free category
$X$ and every small category $Y$, the natural functor $ \bbD(X
\times Y) \longrightarrow \Fun(X^\op,\bbD(Y))$ is full and
essentially surjective.
\item[-] A derivator $\bbD$ is called {\em triangulated} (or {\em stable}) if it is
pointed and if every global commutative square in $\bbD$ is
cartesian exactly when it is cocartesian. A source of examples is provided by the derivators $\HO(\cM)$ associated to {\em stable} Quillen model categories $\cM$.
\end{enumerate}
Recall from \cite{CN} that given any triangulated derivator $\bbD$ and small category $X$, the category $\bbD(X)$ has a canonical triangulated structure. In particular, the category $\bbD(e)$ is triangulated. Recall also from 
\cite{CT} that any triangulated derivator
$\bbD$ is canonically enriched over spectra, \ie we have a well-defined morphism of derivators
$$ \bbR\Hom(-,-): \bbD^\op \times \bbD \too \HO(\Spt)\,.$$

Finally, given derivators $\bbD$ and $\bbD'$, we denote by $\uHom(\bbD,\bbD')$ the
category of all morphisms of derivators and by $\HomC(\bbD,\bbD')$ the
category of morphisms of derivators which preserve arbitrary homotopy
colimits.

\end{document}